\documentclass[12pt]{article}%
\usepackage{amsmath}
\usepackage[dvips]{graphicx}
\usepackage{amsfonts}
\usepackage{amssymb}%
\begin{document}

\title{\normalsize {\bf ALMOST INVARIANT HALF-SPACES FOR OPERATORS
 ON HILBERT SPACE. II:
OPERATOR MATRICES}}
\author{{\normalsize Il Bong Jung, Eungil Ko and Carl Pearcy\thanks{\textit{%
2000 Mathematics Subject Classification}.  Primary 47A15.} \thanks{\textit{%
Key words and phrases: }Invariant subspace, half-space, almost invariant
half-space, operator matrices.}}}
\date{}
\maketitle

\begin{abstract}
This paper is a sequel to \cite{JKP}. In that paper we transferred the
discussions in \cite{APTT} and \cite{PT} concerning almost invariant
half-spaces for operators on complex Banach spaces to the context of
operators on Hilbert space, and we gave easier proofs of the main results in
\cite{APTT} and \cite{PT}. In the present paper we discuss consequences of
the above-mentioned results for the matricial structure of operators on
Hilbert space.
\end{abstract}

\bigskip

\bigskip
\textbf{1. Introduction.} This paper is a sequel to \cite{JKP} and, for
brevity and to avoid repetition, we shall assume that the reader is familiar
with the contents of \cite{JKP}, including the notation, definitions,
results, and proofs.

We begin with a few introductory remarks so as to give sufficient credit to
the authors of \cite{APTT} and \cite{PT}. In those papers the authors
originated, in the context of infinite dimensional, complex Banach spaces, a
completely new and remarkable construction that allowed them to prove that
all operators on many such spaces have an almost invariant half-space of
defect $0$ or $1$. This is as close as anyone has gotten to resolving the
general invariant subspace problem for operators on a Hilbert space, and the
authors of \cite{APTT} and \cite{PT} (and several other authors of
subsequent articles; see the references) have surely advanced operator
theory.

In \cite{JKP} the present authors changed the context of this developing
theory to operators on Hilbert space, where there is more structure that
enabled us to give simpler proofs of the main results of \cite{APTT} and
\cite{PT} (cf. \cite[Theorems 1.4 and 2.1]{JKP}). We also obtained the new
result \cite[Theorem 3.1]{JKP}.

One of the main reasons for that effort is that those new proofs allow us,
in this paper, to obtain some new matricial \textquotedblleft standard
forms\textquotedblright\ for operators on Hilbert space. These results
generalize earlier theorems in \cite{BP1} and \cite{BP2} and may be useful
in several areas of operator theory.

Throughout this article, as in \cite{JKP}, $\mathcal{H}$ will always denote
a separable, infinite dimensional, complex Hilbert space, $\mathcal{L(H)}$
the algebra of all bounded linear operators on $\mathcal{H},$ and $\mathbb{F}
$ the ideal in $\mathcal{L(H)}$ of finite-rank operators. We shall need
several subsets of $\sigma (T)$, the spectrum of $T$, including $\sigma
_{e}(T)$, the essential (Calkin) spectrum, $\partial \sigma _{e}(T)$, the
boundary of the essential spectrum, $\sigma _{lre}(T)$, the intersection of
the left and right essential spectra, and $\sigma _{p}(T)$, the point
spectrum of $T$. We write, as usual, $\mathbb{C}1_{\mathcal{H}}$ for the set
of all scalar multiples of the identity operator $1_{\mathcal{H}},$ and we
remind the reader that operators $A$ and $B$ in $\mathcal{L(H)}$ are called
\textit{equivalent} if there exist invertible operators $S_{1}$ and $S_{2}$
in $\mathcal{L(H)}$ such that $S_{1}AS_{2}=B$. This equivalence relation was
studied by K\"{o}the\cite{Ko1},\cite{Ko2}, and by L. Williams in his thesis%
\cite{Wil} (among others), and in certain cases a complete set of invariants
for this relation are known.

We note that in what follows, when an operator $A:\mathcal{H}_{1}\rightarrow
\mathcal{H}_{2}$, with dim $\mathcal{H}_{1}=$ dim $\mathcal{H}_{2}=\aleph
_{0}$, is said to belong to an ideal $\mathcal{J\subset L(H}_{1}\mathcal{)}$%
, this is to be interpreted to mean that there exists a Hilbert space
isomorphism $\varphi $ of $\mathcal{H}_{2}$ onto $\mathcal{H}_{1}$ such that
when $\mathcal{H}_{2}$ is identified with $\mathcal{H}_{1}$ via this
isomorphisms, then $A\in \mathcal{J}(\mathcal{H}_{1})$.

\bigskip

\textbf{2. The new results. }Our first result, which is new and generalizes
old results in \cite{BP1}, is the following:

\medskip

\textbf{Theorem 2.1.} \textit{Let }$T^{\prime }$\textit{\ be arbitrary in }$%
\mathcal{L(H)}$,\textit{\ let} $\mathcal{J(H)}$ \textit{be any proper ideal
in }$\mathcal{L(H)}$\textit{\ that is strictly larger than the ideal }$%
\mathbb{F}$\textit{, and let }$\varepsilon >0$\textit{\ be given. Then there
exist a scalar }$\beta =\beta (T^{\prime })$,\textit{\ an invertible
operator }$S$\textit{\ in }$\mathcal{L(H)}$, \textit{a half}-\textit{space} $%
\mathcal{M}$ \textit{of} $\mathcal{H}$\textit{, and} \textit{operators }$%
T_{1,1}$ \textit{and }$R$\textit{\ in }$\mathcal{J}$ \textit{with }$R$%
\textit{\ of rank }$0$\textit{\ or }$1$\textit{,} \textit{such that the
matrix of }$S^{-1}(T^{\prime }-\beta 1_{\mathcal{H}})S$ \textit{with respect
to the decomposition $\mathcal{H=M\oplus M}^{\perp }$\ has the form}%
\begin{equation}
\left(
\begin{array}{cc}
T_{1,1} & T_{1,2} \\
R & T_{2,2}%
\end{array}%
\right)  \tag{2.1}
\end{equation}%
\textit{where }

a$)$ $T_{1,1}$\textit{\ is a diagonalizable normal operator with }$0\in
\partial \sigma _{e}(T_{1,1})$ \textit{and }$\Vert T_{1,1}\Vert <\varepsilon
$,

b$)$ $\left\Vert R\right\Vert <\varepsilon $,

c$)$\textit{\ if }$\mathcal{J(M)}$\textit{\ is taken to be the trace-class%
\textit{\ and that norm is denoted by}} $|||~~|||$\textit{$_{\mathcal{J}}$,%
\textit{\ then }}$|||T_{1,1}|||_{\mathcal{J}}<\varepsilon $\textit{\ and} $%
|||R|||_{\mathcal{J}}<\varepsilon $, \textit{and}

d$)$\textit{\ if one gives up the requirement on the smallness of the norms
of }$T_{1,1}$\textit{\ and }$R$,\textit{\ then the\ operator }$T_{1,2}$%
\textit{\ appearing in }$(2.1)$ \textit{can be replaced by any operator
equivalent to it in the sense defined above. In particular, if} dim ker$%
(T_{1,2})=$ dim $\ker (T_{1,2}^{\ast })$, \textit{then }$T_{1,2}$\textit{\
can be replaced by a positive semidefinite diagonalizable operator. If, in
addition, }$T_{1,2}$\textit{\ has closed range, this diagonalizable operator
can be taken to be a projection.}

\smallskip

\textit{Proof}. This proof closely resembles those of \cite[Theorems 1.4 and
2.1]{JKP} but with some additional considerations. As in \cite{JKP}, for
clarity we separate the different parts of the proof into cases.

\smallskip

Case 1. \textit{There exists a scalar }$\beta =\beta (T^{\prime })$\textit{\
in} $\partial \sigma _{e}(T^{\prime })$\textit{\ such that }$\beta \not\in
\sigma _{p}(T^{\prime }).$

\smallskip

In the proof we shall write $T$ for $T^{\prime }-\beta 1_{\mathcal{H}}$, and
thus we know that $0\in \partial \sigma _{e}(T)$\ and $0\not\in \sigma
_{p}(T)$. Let now $\mathcal{J}(\mathcal{H})$ be any proper ideal in $%
\mathcal{L(H)}$ that strictly contains the ideal $\mathbb{F}$, and let $%
\varepsilon >0$ be given. (Those readers wishing to refamiliarize themselves
with the general theory of ideals in $\mathcal{L(H)}$ might consult \cite%
{BPS1}.) Let also $\mathcal{S}_{\mathcal{J}}$ be the ideal set associated
with $\mathcal{J(H)}$, and recall that $\partial \sigma (T)\setminus
\partial \sigma _{e}(T)$ consists of at most a countable number of isolated
eigenvalues of $T$ whose derived set is a subset of $\partial \sigma _{e}(T)$%
. Thus we may choose a sequence of distinct points $\{\lambda _{n}\}_{n\in
\mathbb{N}}\subset \mathbb{C}\setminus \sigma (T)$ that converges to $0$
sufficiently fast that $\left\{ \left\vert \lambda _{n}\right\vert \right\}
_{n\in \mathbb{N}}\in \mathcal{S}_{\mathcal{J}}$. Note that the sequence $%
\{\Vert \left( T-\lambda _{n}1_{\mathcal{H}}\right) ^{-1}\Vert \}_{n\in
\mathbb{N}}$ converges to $+\infty $, and the uniform boundedness theorem
provides a unit vector $e$ and a subsequence $\{\lambda _{n_{k}}\}_{k\in
\mathbb{N}}$ such that
\begin{equation*}
\lim_{k\rightarrow \infty }\Vert (T-\lambda _{n_{k}}1_{\mathcal{H}%
})^{-1}e\Vert =+\infty .
\end{equation*}%
Note also that by one of the defining properties of an ideal set (cf. \cite%
{BP1}), $\left\{ \left\vert \lambda _{n_{k}}\right\vert \right\} _{k\in
\mathbb{N}}\in \mathcal{S}_{\mathcal{J}}$. To conserve notation, we rename
this last sequence $\{\lambda _{n}\}_{n\in \mathbb{N}}$, so that
\begin{equation*}
\lim_{n\rightarrow \infty }\Vert (T-\lambda _{n}1_{\mathcal{H}})^{-1}e\Vert
=+\infty .
\end{equation*}%
We next define sequences $\{\alpha _{n}\}_{n\in \mathbb{N}}\subset \mathbb{R}%
_{+}$ and $\{h_{n}\}_{n\in \mathbb{N}}\subset \mathcal{H}$ by the formulas%
\begin{equation*}
\alpha _{n}=\Vert (T-\lambda _{n}1_{\mathcal{H}})^{-1}e\Vert ^{-1},~\ \
h_{n}=\alpha _{n}(T-\lambda _{n}1_{\mathcal{H}})^{-1}e,~\ \ n\in \mathbb{N},
\end{equation*}%
and note that $\{\alpha _{n}\}\rightarrow 0$ and the unit vectors $h_{n}$
satisfy%
\begin{equation}
Th_{n}=\lambda _{n}h_{n}+\alpha _{n}e,~\ \ \ n\in \mathbb{N}.  \tag{2.2}
\end{equation}

As in the proof of \cite[Theorem 2.1]{JKP}, we extract a subsequence $%
\{h_{n_{k}}\}_{k\in \mathbb{N}}$ of $\{h_{n}\}_{n\in \mathbb{N}}$ that
converges weakly to some vector $h^{\prime }\in \mathcal{H}$ and has the
additional properties that $\sum_{k\in \mathbb{N}}\alpha
_{n_{k}}^{2}<+\infty $ and $\{\alpha _{n_{k}}\}_{k\in \mathbb{N}}\in
\mathcal{S}_{\mathcal{J}}\,.$ We also replace $h_{n_{1}}$ by $e$, so that
henceforth $h_{n_{1}}=e$ (but (2.2) is not satisfied for $n=1$). Equation
(2.2) plus the fact that $0\not\in \sigma _{p}(T)$ shows that $h^{\prime }=0$%
, and we harmlessly rename the sequence $\{h_{n_{k}}\}_{k\in \mathbb{N}}$ as
$\{l_{n}\}_{n\in \mathbb{N}}$, so $l_{1}=e$ and $\{l_{n}\}_{n\in \mathbb{N}}$
converges weakly to $0$. We next use this last fact to construct by
induction a subsequence $\{l_{n_{k}}\}_{k\in \mathbb{N}}$ of $%
\{l_{n}\}_{n\in \mathbb{N}}$ such that $l_{n_{1}}=l_{1}=e$ and
\begin{equation}
\left\vert \left\langle l_{n_{j}},l_{n_{m}}\right\rangle \right\vert
<1/4^{j+m},~\ \ \ \ j,m\in \mathbb{N},~\ \ j\not=m.  \tag{2.3}
\end{equation}%
We show now that the sequence $\{l_{n_{k}}\}_{k\in \mathbb{N}}$ thus
constructed is linearly independent. Suppose, to the contrary, that some
equation of the form%
\begin{equation*}
l:=\sum_{p=1}^{N}\gamma _{p}l_{n_{j(p)}}=0,
\end{equation*}%
where all $\gamma _{p}$ are nonzero and $\sum_{p=1}^{N}\left\vert \gamma
_{p}\right\vert ^{2}=1$. Upon taking the inner product of $l$ with itself,
one obtains%
\begin{equation*}
1=\sum_{p=1}^{N}\left\vert \gamma _{p}\right\vert ^{2}=-\sum_{\substack{ %
p,q=1 \\ p\not=q}}^{N}2\rm{Re}\left\langle
l_{n_{j(p)}},l_{n_{j(q)}}\right\rangle ,
\end{equation*}%
which, after a bit of arithmetic, clearly contradicts (2.3) and gives the
linear independence. We can now define an orthonormal sequence $%
\{f_{k}\}_{k\in \mathbb{N}}$ by induction (or by applying the Gram-Schmidt
procedure), setting $f_{1}=l_{n_{1}}(=e)$ and, in general, assuming that
orthonormal vectors $\left\{ f_{1},...,f_{j}\right\} $ have been defined
such that $\vee \left\{ f_{1},...,f_{j}\right\} =\vee \left\{
l_{n_{1}},...,l_{n_{j}}\right\} $, we define the unit vector $f_{j+1}\in
\vee \left\{ l_{n_{1}},...,l_{n_{j+1}}\right\} $ to be orthogonal to $\vee
\left\{ f_{1},...,f_{j}\right\} $. Also set $\mathcal{M}^{\prime }=\vee
_{k\in \mathbb{N}}\{f_{k}\}=\vee _{k\in \mathbb{N}}\{l_{n_{k}}\}$, and
define a linear transformation $S:\mathcal{M}^{\prime }\rightarrow \mathcal{M%
}^{\prime }\ $by$\ S(\sum_{k\in \mathbb{N}}\beta _{k}f_{k})=\sum_{k\in
\mathbb{N}}\beta _{k}l_{n_{k}}$ for every finitely nonzero sequence $\{\beta
_{k}\}_{k\in \mathbb{N}}$. The inequality (2.3), together with an easy
calculation, shows that whenever $\sum_{k\in \mathbb{N}}$ $\left\vert \beta
_{k}\right\vert ^{2}=1$ we have
\begin{equation*}
43/45\leq \Vert \sum_{k\in \mathbb{N}}\beta _{k}l_{n_{k}}\Vert ^{2}\leq
47/45,
\end{equation*}%
and thus $S$ extends by continuity to a bounded, invertible operator in $%
\mathcal{L(M}^{\prime })$ such that $Sf_{k}=l_{n_{k}},$ $k\in \mathbb{N}$.

Furthermore we define $S$ on $(\mathcal{M}^{\prime })^{\perp }$ to be $1_{(%
\mathcal{M}^{\prime })^{\perp }},$ so that henceforth $S$ is an invertible
operator in $\mathcal{L(H)}$. Now let us examine the action of the operator $%
S^{-1}(T^{\prime }-\beta 1_{\mathcal{H}})S$ on the orthonormal basis $%
\{f_{n}\}_{n\in \mathbb{N}}$ for $\mathcal{M}^{\prime }$. We have from (2.2)
that
\[
\begin{aligned}
S^{-1}(T^{\prime }-\beta 1_{\mathcal{H}})S|_{\mathcal{M}^{\prime }}(f_{k})
&=S^{-1}T(l_{n_{k}})  \notag \\
&=S^{-1}(\lambda _{n_{k}}l_{n_{k}}+\alpha _{n_{k}}e)   \\
&=\lambda _{n_{k}}f_{k}+\alpha _{n_{k}}f_{1},~\ \ \ \ \ \ k\in \mathbb{N}%
\setminus \{1\}.
\end{aligned}\eqno (2.4)
\]
Moreover, defining $P_{\mathcal{M}^{\prime }}$ to be the projection of $%
\mathcal{H}$ onto $\mathcal{M}^{\prime }$, we have
\begin{equation*}
S^{-1}P_{\mathcal{M}^{\prime }}(T^{\prime }-\beta 1_{\mathcal{H}})S|_{%
\mathcal{M}^{\prime }}(f_{1})=S^{-1}P_{\mathcal{M}^{\prime }}Te.
\end{equation*}%
The presence of $P_{\mathcal{M}^{\prime }}$ in this equation is necessary
because $Te$ may not belong to $\mathcal{M}^{\prime }$. If $Te\in \mathcal{M}%
^{\prime }$, then $\mathcal{M}^{\prime }$ is an invariant subspace for $T$
and we will see shortly below that $\mathcal{M}^{\prime }$ contains an
invariant half-space for $T$. Thus if $Te\in \mathcal{M}^{\prime }$, then $%
R=0$ and the proof is complete in that case. Thus we may assume henceforth
that $Te\not\in \mathcal{M}^{\prime }$. It follows easily from (2.4) that
the matrix for $S^{-1}P_{\mathcal{M}^{\prime }}TS|_{\mathcal{M}^{\prime }}$
with respect to the orthonormal basis $\{f_{k}\}_{k\in \mathbb{N}}$ for $%
\mathcal{M}^{\prime }$ is
\begin{equation}
\left(
\begin{array}{ccccc}
\beta _{1} & \alpha _{n_{2}} & \alpha _{n_{3}} & \alpha _{n_{4}} & \cdots
\\
\beta _{2} & \lambda _{n_{2}} & 0 & 0 & \cdots  \\
\beta _{3} & 0 & \lambda _{n_{3}} & 0 & \cdots  \\
\beta _{4} & 0 & 0 & \lambda _{n_{4}} &  \\
\vdots  & \vdots  & \vdots  &  & \ddots
\end{array}%
\right) ,  \tag{2.5}
\end{equation}%
where $P_{\mathcal{M}^{\prime }}Te=\sum_{j\in \mathbb{N}}\beta _{j}l_{n_{j}}$%
.

We are finally ready to define the half-space $\mathcal{M}$ so that the
matrix for $S^{-1}(T^{\prime }-\beta 1_{\mathcal{H}})S$ relative to the
decomposition $\mathcal{H}=\mathcal{M\oplus M}^{\perp }$ has the properties
called for in the statement of the theorem. We already know that both
sequences $\{|\lambda _{n_{k}}|\}_{k\in \mathbb{N}}$ and $\{\alpha
_{n_{k}}\}_{k\in \mathbb{N}}$ belong to $\mathcal{S}_{\mathcal{J}}$ so to
show that $T_{1,1}$ and $R$ satisfy the required norm inequalities, it
suffices to ignore the first few $f_{k}$ $(k>1)$ when defining $\mathcal{M}$%
. Thus we choose an appropriate positive integer $K$ and define $\mathcal{M}%
=\vee _{k\geq K}f_{2k},$ which ensures that $\mathcal{M}$ is a half-space of
$\mathcal{H}$. A careful consideration of the matrices (2.5) and (2.1) now
shows that this completes the proof of a), b), and c) of the theorem in Case
1. (Here we are implicitly using the fact that $\mathcal{M}^{\perp }$ may be
identified with a second copy of $\mathcal{M}$ to make sense of statements
like $R\in \mathcal{J}$.) To show that d) can be satisfied if the
requirement that the norms of $T_{1,1}$ and $R$ be small is dropped, let $%
S_{1}\in \mathcal{L(M)}$ and $S_{2}\in \mathcal{L(M}^{\perp }\mathcal{)}$ be
arbitrary invertible operators and take note of the matricial calculation%
\begin{equation}
\left(
\begin{array}{cc}
S_{1} & 0 \\
0 & S_{2}^{-1}%
\end{array}%
\right) \left(
\begin{array}{cc}
T_{1,1} & T_{1,2} \\
R & T_{2,2}%
\end{array}%
\right) \left(
\begin{array}{cc}
S_{1}^{-1} & 0 \\
0 & S_{2}%
\end{array}%
\right) =\left(
\begin{array}{cc}
S_{1}T_{1,1}S_{1}^{-1} & S_{1}T_{1,2}S_{2} \\
S_{2}^{-1}RS_{1}^{-1} & S_{2}^{-1}T_{2,2}S_{2}%
\end{array}%
\right) .  \tag{2.6}
\end{equation}%
Equation (2.6) now shows that $T_{1,2}$ may be replaced by any operator
equivalent to it. In particular, that the remaining statements in d) can be
satisfied is immediate from the work of K\"{o}the\cite{Ko1},\cite{Ko2} (see
also \cite{Wil}). The proof of Case 1 is now complete.

\smallskip

Case 2. $\partial \sigma (T^{\prime })$\textit{\ is an infinite set and }$%
\partial \sigma _{e}(T^{\prime })\subset \sigma _{p}(T^{\prime }).$

\smallskip

Recall from the Fredholm theory that every point of $\partial \sigma
(T^{\prime })\setminus \partial \sigma _{e}(T^{\prime })$ is an isolated
eigenvalue of $T^{\prime }$ of finite multiplicity and that every
accumulation point of this set belongs to $\partial \sigma _{e}(T^{\prime })$%
. Thus there exist $\beta \in \partial \sigma _{e}(T^{\prime }),$ a sequence
$\{\lambda _{n}\}$ of nonzero eigenvalues of $T^{\prime }-\beta 1_{\mathcal{H%
}}$ converging to $0$, and a corresponding sequence $\{v_{n}\}$ of
eigenvectors of $T^{\prime }-\beta 1_{\mathcal{H}}$. We define $\mathcal{M}%
=\vee _{n\in \mathbb{N}}\{v_{n}\}$, note that $\mathcal{M}$ is an infinite
dimensional invariant subspace of $T^{\prime }-\beta 1_{\mathcal{H}}$, and
that $((T^{\prime }-\beta 1_{\mathcal{H}})\mathcal{M})^{-}=\mathcal{M}$.
Thus $0\in \partial \sigma _{e}(((T^{\prime }-\beta 1_{\mathcal{H}})|_{%
\mathcal{M}})^{\ast })$ but $0\not\in \sigma _{p}(((T^{\prime }-\beta 1_{%
\mathcal{H}})|_{\mathcal{M}})^{\ast })$. In other words, $\left( (T^{\prime
}-\beta 1_{\mathcal{H}})|_{\mathcal{M}}\right) ^{\ast }$ satisfies the
hypotheses of Case 1 of this theorem, which has already been proved. Thus $%
\left( (T^{\prime }-\beta 1_{\mathcal{H}})|_{\mathcal{M}}\right) ^{\ast }$
satisfies all of the conclusions of the theorem. And it is easy to see that
those conclusions are preserved by the mappings $A^{\ast }\rightarrow A$ and
$A|_{\mathcal{M}}\rightarrow A,$ $A\in \mathcal{L(H)}$. Thus $T^{\prime
}-\beta 1_{\mathcal{H}}$ satisfies the conclusions of the theorem, which
completes the proof of Case 2.

\smallskip

Case 3. $\partial \sigma (T^{\prime })$\textit{\ is a finite set. }

\smallskip

This obviously implies that $\sigma (T^{\prime })$ is finite, and thus there
exists $\beta \in \partial \sigma _{e}(T^{\prime })$ that is an isolated
point of $\sigma (T^{\prime })$. One quickly reduces, as in \cite{JKP}, to
the case in which $T^{\prime }-\beta 1_{\mathcal{H}}$ is a quasinilpotent
operator acting on an infinite dimensional space. If $T^{\prime }-\beta 1_{%
\mathcal{H}}$ is nilpotent, the result is trivially true since ker$%
(T^{\prime }-\beta 1_{\mathcal{H}})$ is an infinite dimensional invariant
subspace for $T^{\prime }-\beta 1_{\mathcal{H}}$ which, if necessary, can be
\textquotedblleft split into halves\textquotedblright\ to ensure that $%
\mathcal{M}\subset $ ker$(T^{\prime }-\beta 1_{\mathcal{H}})$ is a
half-space with $(T^{\prime }-\beta 1_{\mathcal{H}})|_{\mathcal{M}}=0$.
Setting both $T_{11}$ and $R$ equal $0$ completes the argument when $%
T^{\prime }-\beta 1_{\mathcal{H}}$ is nilpotent. Thus henceforth we suppose
that $T^{\prime }-\beta 1_{\mathcal{H}}$ is quasinilpotent but not nilpotent
and that $0\in \sigma _{p}(T^{\prime }-\beta 1_{\mathcal{H}})$. If either ker%
$(T^{\prime }-\beta 1_{\mathcal{H}})$ or ker$(T^{\prime ^{\ast }}-\bar{\beta}%
1_{\mathcal{H}})$ is infinite dimensional, one can, by applying \cite[Remark
1.3]{JKP} if necessary (or, equivalently, using the fact, mentioned above,
that the map $A^{\ast }\rightarrow A$ preserves the conclusions of the
present theorem), construct the matrix (2.1) for $T^{\prime }-\beta 1_{%
\mathcal{H}}$ with $T_{1,1}=R=0$. Moreover, if $0\not\in \sigma
_{p}(T^{\prime }-\beta 1_{\mathcal{H}})\cap \sigma _{p}((T^{\prime }-\beta
1_{\mathcal{H}})^{\ast })$, then, again using the just mentioned trick, we
are exactly in Case 1 with perhaps $(T^{\prime }-\beta 1_{\mathcal{H}%
})^{\ast }$ replacing $T^{\prime }-\beta 1_{\mathcal{H}}$. Thus, to complete
the proof of the theorem, it suffices to dispose of the case in which both
ker$(T^{\prime }-\beta 1_{\mathcal{H}})$ and ker$((T^{\prime }-\beta 1_{%
\mathcal{H}})^{\ast })$ are nonzero and finite dimensional. Observe next
that it suffices to exhibit an infinite dimensional subspace $\mathcal{N}%
\subset \mathcal{H}$, invariant for $(T^{\prime }-\beta 1_{\mathcal{H}%
})^{\ast }$, such that $(T^{\prime }-\beta 1_{\mathcal{H}})^{\ast }|_{%
\mathcal{N}}$ has range dense in $\mathcal{N}$. This is because $(T^{\prime
}-\beta 1_{\mathcal{H}})^{\ast }$ is quasinilpotent and thus so is $%
(T^{\prime }-\beta 1_{\mathcal{H}})^{\ast }|_{\mathcal{N}}$, and such an $%
\mathcal{N}$ would have the property that $\sigma _{p}((T^{\prime }-\beta 1_{%
\mathcal{H}})^{\ast }|_{\mathcal{N}})^{\ast }=\emptyset $, which is exactly
Case 1 again. Finally, the argument that produces a subspace $\mathcal{N}$
with the desired properties is written in great detail in \cite[Proof of
Theorem 1.4, Case III]{JKP} so we do not repeat it here. This concludes the
proof of Theorem 2.1. \ $\square $

\medskip

\textbf{Remark 2.2.} The presence of the similarity $S^{-1}[~~]S$ in the
above theorem is unpleasant, but it has two functions, neither of which can
seemingly be avoided. It is needed to orthogonalize the \textquotedblleft
eigenvectors\textquotedblright\ $l_{n_{k}}$ and to ensure that the vector $e$
is orthogonal to those \textquotedblleft eigenvectors\textquotedblright .
One may eliminate the similarity $S^{-1}[~~]S$ from Theorem 1.1 at the cost
of replacing the orthogonal direct sum $\mathcal{M}\oplus \mathcal{M}^{\perp
}$ by an \textit{algebraic direct sum }$\mathcal{M}\overset{\mathbf{\bullet }%
}{+}\mathcal{N}$. Recall that to say that $\mathcal{H}=\mathcal{M}\overset{%
\mathbf{\bullet }}{+}\mathcal{N}$ means that $\mathcal{H}=\mathcal{M}+%
\mathcal{N}$ and $\mathcal{M}\cap \mathcal{N}=(0).$

\medskip

Relative to such a decomposition of $\mathcal{H}$, Theorem 2.1 becomes:

\medskip

\textbf{Corollary 2.3.} \textit{With the same hypotheses as in Theorem 2.1,
there exist half-spaces }$\mathcal{M}$ \textit{and }$\mathcal{N}$ \textit{of}
$\mathcal{H}$ \textit{with} $\mathcal{H}=\mathcal{M\overset{\mathbf{\bullet }%
}{+}N}$ \textit{such that the matrix of} $(T^{\prime }-\beta 1_{\mathcal{H}%
}) $ \textit{with respect to this decomposition has the form}
\begin{equation}
\left(
\begin{array}{cc}
\widehat{T}_{1,1} & \widehat{T}_{1,2} \\
\widehat{R} & \widehat{T}_{2,2}%
\end{array}%
\right) ,  \tag{2.7}
\end{equation}%
\textit{where}

a$)$ $\widehat{T}_{1,1}\in \mathcal{J(M)}$\textit{\ and }$\widehat{T}_{1,1}$
\textit{is similar to a diagonalizable normal operator,}

b$)$ $\widehat{R}:\mathcal{M\rightarrow N}$ \textit{is a rank-one operator,}

c$)$ \textit{the operator }$\widehat{T}_{1,2}:\mathcal{N\rightarrow M}$%
\textit{\ may be replaced by any operator equivalent to it.}

\medskip

\textbf{Remark 2.4.} With a more delicate argument, the norms of \textit{\ }$%
\widehat{T}_{1,1}$ and $\widehat{R}$ in Corollary 2.3 can also be made small.

\medskip

The fact that every nonscalar operator $T^{\prime }$ in $\mathcal{L(H)}$ is
associated with a matrix of the form (2.1) raises some interesting, and
possibly important, problems. Here we pose some of them.

\medskip

\textbf{Problem 2.5.} Suppose $T$ is a non-normal hyponormal operator in $%
\mathcal{L(H)}$. Then, via Theorem 2.1 above, there exists a half-space $%
\mathcal{M}$ of $\mathcal{H}$ and an invertible $S\in \mathcal{L(H)}$ such
that the matrix of $S^{-1}TS$ with respect to the decomposition $\mathcal{H}=%
\mathcal{M\oplus M}^{\perp }$ has the form%
\begin{equation*}
T=\left(
\begin{array}{cc}
T_{11} & T_{12} \\
R & T_{22}%
\end{array}%
\right) ,
\end{equation*}%
where $R$ has rank at most $1$. Can this matrix for $S^{-1}TS$ be used to
help solve the (open) invariant subspace problem for hyponormal operators?

\medskip

\textbf{Problem 2.6.} Is it possible that, with some clever argument, one
could show that the operator $T_{1,2}$ in (2.1) can frequently be taken to
be of finite rank? Said differently, can many operators $T$ in $\mathcal{L(H)%
}$ be written as a sum of a finite-rank operator and an operator $\widetilde{%
T}$ with the property that there exist invariant half-spaces $\mathcal{M}$
and $\mathcal{N}$ for $\widetilde{T}$ such that $\mathcal{M}\overset{\bullet
}{+}\mathcal{N}=\mathcal{H}$.

\bigskip

\textbf{3. Adjoint considerations.} In this section we exhibit the
modifications that can be made to the matrix in (2.1) if, after applying
Theorem 2.1 to an operator $T^{\prime }$ in $\mathcal{L(H)}$, we then apply
it to $T_{2,2}^{\ast }$ in (2.1).

\medskip

\textbf{Theorem 3.1.} \textit{Let }$T^{\prime }$\textit{\ be an arbitrary
operator in} $\mathcal{L(H)}$, \textit{let }$\varepsilon >0$\textit{\ be
arbitrary, and let} $\mathcal{J(H)}$ \textit{be a proper ideal in }$\mathcal{%
L(H)}$ \textit{strictly larger than} $\mathbb{F}$. \textit{Then there exist
a scalar }$\beta $\textit{, an invertible operator }$S$ \textit{in $\mathcal{%
L(H)}$}, \textit{and infinite dimensional orthogonal subspaces} $\mathcal{M}%
_{i}$, $i=1,2,3,$ \textit{in} $\mathcal{H}$ \textit{such that} $\mathcal{M}%
_{1}\oplus \mathcal{M}_{2}\oplus \mathcal{M}_{3}=\mathcal{H}$ \textit{and
such that the matrix for} $S^{-1}(T^{\prime }-\beta 1_{\mathcal{H}})S$
\textit{relative to this decomposition has the form}%
\begin{equation}
\left(
\begin{array}{ccc}
T_{1,1} & T_{1,2} & T_{1,3} \\
R_{2,1} & T_{2,2} & T_{2,3} \\
R_{3,1} & R_{3,2} & T_{3,3}%
\end{array}%
\right) ,  \tag{3.1}
\end{equation}%
\textit{where}

a$)$ $T_{1,1}\in \mathcal{J(M}_{1})$, $T_{3,3}\in \mathcal{J(M}_{3})$,
\textit{and }$\left\Vert T_{1,1}\right\Vert <\varepsilon $, $\left\Vert
T_{3,3}\right\Vert <\varepsilon $,

b$)$ $R_{2,1},R_{3,1}$\textit{, and }$R_{3,2}$ \textit{are operators of rank
at most one, and}

c$)$ $T_{1,1}$\textit{\ and }$T_{3,3}$\textit{\ are diagonalizable normal
operators. }

\smallskip

\textit{Proof.} With $T^{\prime }$, $\varepsilon $, and $\mathcal{J}$ as in
the hypothesis, choose $\beta \in \partial \sigma _{e}(T^{\prime })$ and $%
S\in \mathcal{L(H)}$ invertible to satisfy the conclusions of Theorem 2.1.
Thus there exists a decomposition of $\mathcal{H}$ into half-spaces $%
\mathcal{M}_{1}$ and $\mathcal{M}_{2}$ with $\mathcal{H}=\mathcal{M}%
_{1}\oplus \mathcal{M}_{2}$ such that relative to this decomposition, the
matrix for $S^{-1}TS$ is as in (2.1) where $T:=T^{\prime }-\beta 1_{\mathcal{%
H}}$. Let $\mathcal{E}=\{e_{n}\}_{n\in \mathbb{N}}$ be the orthonormal basis
for $\mathcal{M}_{1}$ relative to which the matrix for $T_{1,1}$ is written,
and recall that there exists a sequence $\{\lambda _{n}\}_{n\in \mathbb{N}%
}\subset \mathbb{C}$ such that $T_{11}e_{n}=\lambda _{n}e_{n},$ $n\in
\mathbb{N}$ (all $\lambda _{n}$ may be $0)$. This shows that $0\in \sigma
_{le}(T_{11})$, and since $0\in \partial \sigma _{e}(T)$, we have $0\in
\partial \sigma _{e}(T_{1,1})$. Define $\mathcal{N}_{1}=\vee _{n\in \mathbb{N%
}}\{e_{2n}\}$ and $\mathcal{N}_{2}=\vee _{n\in \mathbb{N}}\{e_{2n-1}\}$, so $%
\mathcal{N}_{1}\oplus \mathcal{N}_{2}=\mathcal{M}_{1}$. Note that $0\in
\partial \sigma _{e}(T_{1,1}|_{\mathcal{N}_{i}}),$ $i=1,2$. Now redecompose $%
\mathcal{H}$ as $\mathcal{H}=\mathcal{N}_{1}\oplus (\mathcal{N}_{2}\oplus
\mathcal{M}_{2})$ and observe that the matrix $(\widetilde{T}_{i,j})$ for $%
S^{-1}TS$ relative to this decomposition has all the properties a)-d) in
Theorem 2.1 and $\widetilde{T}_{1,1}=T_{1,1}|_{\mathcal{N}_{1}}$. What we
have gained by this new decomposition of $\mathcal{H}$ is that now $%
\widetilde{T}_{2,2}\in \mathcal{L(N}_{2}\mathcal{\oplus M}_{2}\mathcal{)}$
and $0\in \partial \sigma _{e}(\widetilde{T}_{2,2})$. It obviously follows
that also $0\in \partial \sigma _{e}((\widetilde{T}_{2,2})^{\ast })$, and we
now apply Theorem 2.1 to the operator $(\widetilde{T}_{2,2})^{\ast }$ with $%
\beta =0$. Then $\mathcal{N}_{2}\oplus \mathcal{M}_{2}$ has a decomposition
as $\mathcal{K}_{1}\oplus \mathcal{K}_{2}$, where the matrix of $(\widetilde{%
T}_{2,2})^{\ast }$ relative to this decomposition of $\mathcal{K}_{1}\oplus
\mathcal{K}_{2}$ has the form%
\begin{equation}
\left(
\begin{array}{cc}
\widehat{T}_{1,1} & \widehat{T}_{1,2} \\
\widehat{R} & \widehat{T}_{2,2}%
\end{array}%
\right)   \tag{3.2}
\end{equation}%
and the entries of this matrix satisfy conditions a)-d) of Theorem 2.1. Thus
there exists another invertible operator $S_{1}$ in $\mathcal{L(H)}$ such
that the matrix of $S_{1}^{-1}(T^{\prime }-\beta 1_{\mathcal{H}})S_{1}$
relative to the decomposition $\mathcal{H}=\mathcal{N}_{1}\oplus \mathcal{K}%
_{2}\oplus \mathcal{K}_{1}$ (note the interchange of $\mathcal{K}_{1}$ and $%
\mathcal{K}_{2}$) is
\begin{equation}
\left(
\begin{array}{ccc}
T_{1,1}|_{\mathcal{N}_{1}} & T_{1,2}^{\prime \prime } & T_{1,3}^{\prime
\prime } \\
R_{2,1} & (\widehat{T}_{2,2})^{\ast } & (\widehat{T}_{1,2})^{\ast } \\
R_{3,1} & (\widehat{R})^{\ast } & (\widehat{T}_{1,1})^{\ast }%
\end{array}%
\right)   \tag{3.3}
\end{equation}%
and a little checking shows that the entries of the matrix for $%
S_{2}^{-1}(T^{\prime }-\beta 1_{\mathcal{H}})S_{2}$ in (3.3) can be made to
have, by careful construction of one more similarity transformation $S_{2}$,
all of the desired properties of the statement of the theorem. In
particular, $T_{1,1}|_{\mathcal{N}_{1}},R_{2,1},R_{3,1},(\widehat{R})^{\ast
},$ and $(\widehat{T}_{1,1})^{\ast }$ all belong to the ideal $\mathcal{J}$
and have operator norms less than $\varepsilon $. Moreover, $T_{1,1}|_{%
\mathcal{N}_{1}}$ and $(\widehat{T}_{1,1})^{\ast }$ are similar to
diagonalizable normal operators. This concludes the proof of the theorem. \ $%
\square $

\medskip

This last corollary is an old theorem of Stampfli\cite{St} with a different
proof that we include as an application of Theorem 3.1 (cf. \cite{BP2}).

\medskip

\textbf{Corollary 3.2.} \textit{Let }$T^{\prime }$\textit{\ be an arbitrary
operator in} $\mathcal{L(H)}$ \textit{and let }$\Delta _{T^{\prime }}$%
\textit{\ be the operator} $($\textit{i.e., derivation}$)$\textit{\ on} $%
\mathcal{L(H)}$ \textit{defined by }$\Delta _{T^{\prime }}(X)=T^{\prime
}X-XT^{\prime }$\textit{\ for }$X$\textit{\ in }$\mathcal{L(H)}.$ \textit{%
Then the range of }$\Delta _{T^{\prime }}$\textit{\ contains no ideal in $%
\mathcal{L(H)}$ strictly larger than} $\mathbb{F}$.

\smallskip

\textit{Proof}. Let $\mathcal{J}_{1}\mathcal{(H)}$ be any ideal in $\mathcal{%
L(H)}$ that properly contains $\mathbb{F}$. One knows from \cite{BPS1} that
there exists an ideal $\mathcal{J}_{2}(\mathcal{H})$ properly containing $%
\mathbb{F}$ such that $\mathcal{J}_{1}(\mathcal{H)}$ properly contains $%
\mathcal{J}_{2}(\mathcal{H)}$. We will now show that the range of $\Delta
_{T^{\prime }}$ does not contain the ideal $\mathcal{J}_{1}(\mathcal{H})$.
With no loss of generality (by applying a unitary equivalence between $%
\mathcal{H}$ and $\mathcal{H}^{(3)}$) we may suppose that $T^{\prime }\in
\mathcal{L(H}^{(3)})$ and (by an application of Theorem 3.1) that there
exist a scalar $\beta $ and an invertible operator $S\in \mathcal{L(H}^{(3)})
$ such that the matrix $M(S^{-1}(T^{\prime }-\beta 1_{\mathcal{H}^{(3)}})S)$
has the form
\begin{equation*}
\left(
\begin{array}{ccc}
T_{1,1} & T_{1,2} & T_{1,3} \\
R_{2,1} & T_{2,2} & T_{2,3} \\
R_{3,1} & R_{3,2} & R_{3,3}%
\end{array}%
\right) ,
\end{equation*}%
where all the entries of this matrix belong to $\mathcal{L(H)}$ and have the
properties a)-c) of Theorem 3.1 (with each of $\mathcal{M}_{1},\mathcal{M}%
_{2},$ and $\mathcal{M}_{3}$ sets equals to $\mathcal{H}$). Now let $X$ be an
arbitrary operator in $\mathcal{L(H}^{(3)})$ having the matrix $%
M(X)=(X_{i,j})_{1\leq i,j\leq 3}$ with entries $X_{i,j}$ in $\mathcal{L(H)}$%
. An elementary calculation shows that
\begin{equation*}
M(\widehat{S})M(X)-M(X)M(\widehat{S})=M(\widehat{S}X-X\widehat{S})
\end{equation*}%
with $\widehat{S}:=S^{-1}(T^{\prime }-\beta 1_{\mathcal{H}^{(3)}})S,$ has
its (3,1) entry in the ideal $\mathcal{J}_{2}(\mathcal{H)\subsetneqq J}_{1}(%
\mathcal{H})$, which clearly shows that the range of $\Delta
_{S^{-1}(T^{\prime }-\beta 1_{\mathcal{H}^{(3)}})S}$ does not contain the
ideal $\mathcal{J}_{1}(\mathcal{H}^{(3)})$. The obvious fact that%
\begin{equation*}
\Delta _{S^{-1}(T^{\prime }-\beta 1_{\mathcal{H}^{(3)}})S}=S^{-1}\Delta
_{T^{\prime }}S
\end{equation*}%
shows that the range of $\Delta _{T^{\prime }}$ cannot contain the ideal $%
\mathcal{J}_{1}(\mathcal{H}^{(3)})$ and completes the proof. \ $\square
$

\medskip

\textbf{Remark 3.3.} We take this opportunity to point out some misprints in
\cite{JKP}. In the statement of Theorem 1.5 the phrase \textquotedblleft $%
\lambda \not\in \sigma _{p}(T)\cap \sigma _{p}(T^{\ast })$%
\textquotedblright\ should read \textquotedblleft $\lambda \not\in \sigma
_{p}(T)$ or $\bar{\lambda}\not\in \sigma _{p}(T^{\ast })$\textquotedblright
. The same misprint occurs in the first line of the proof of Theorem 1.4.

\bigskip

\textbf{Acknowledgement.} The first author was supported by the Basic
Science Research Program through the National Research Foundation of
Korea(NRF) funded by the Ministry of Education (2015R1A2A2A01006072). The
second author was supported by the Basic Science Research Program through
the National Research Foundation of Korea (NRF) funded by the Ministry of
Education (2016R1D1A1B03931937).

\bigskip

\bigskip

\noindent I. B. Jung

\noindent Department of Mathematics, Kyungpook National University, Daegu
41566, Korea

\noindent E-mail: ibjung@knu.ac.kr

\bigskip

\noindent E. Ko

\noindent Department of Mathematics, Ewha Womans University, Seoul 120-750,
Korea

\noindent E-mail: eiko@ewha.ac.kr

\bigskip

\noindent C. Pearcy

\noindent Department of Mathematics, Texas A\&M University, College Station,
TX 77843, USA

\noindent E-mail: pearcy@math.tamu.edu

\end{document}